\documentclass[11pt, oneside]{article}

\usepackage{amsmath, amssymb, amsthm}
\usepackage{enumitem}
\usepackage{capt-of}
\usepackage{ifthen}
\usepackage{geometry}
\usepackage{epsfig}
\usepackage[font=small,labelfont=bf, width=.7\textwidth]{caption}

\newtheorem{theorem}{Theorem}

\newtheorem{claim}{Claim}

\newenvironment{proofc}{{\noindent \textit{Proof of Claim.}}}
{\hfill $\Box$ \vspace*{0.1in}}

\long\def\symbolfootnote[#1]#2{\begingroup\def\thefootnote{\fnsymbol{footnote}}
\footnote[#1]{#2}\endgroup}

\let\oldmarginpar\marginpar
\renewcommand\marginpar[1]{\-\oldmarginpar[\raggedleft\footnotesize #1]%
{\raggedright\footnotesize #1}}

\begin{document}

\title{On a limit of the method of Tashkinov trees for edge-colouring}
\author{
John Asplund\footnote{jasplund@daltonstate.edu; Department of Technology and Mathematics, Dalton State College, Dalton, GA 30720}
\and
Jessica McDonald\footnote{mcdonald@auburn.edu; Department of Mathematics and Statistics,
Auburn University, Auburn, AL 36849}
}
\date{}

\maketitle

\begin{abstract}
The main technique used to edge-colour graphs requiring $\Delta+2$ or more colours is the method of Tashkinov trees. We present a specific limit to this method, in terms of Kempe changes. We also provide a new Tashkinov tree extension.
\end{abstract}

\section{Introduction}

Let $G$ be a loopless graph. The chromatic index of $G$, denoted $\chi'(G)$, is the minimum number $k$ of colours needed to $k$-edge-colour the graph -- that is, to assign colours $\{1, 2, \ldots, k\}$ to the edges of $G$ so that adjacent edges receive different colours. It is obvious that $\chi'(G) \geq \Delta$, where $\Delta$ denotes the maximum degree of $G$. If $G$ is simple then Vizing's Theorem \cite{V2} tells us that $\chi'(G)\in\{\Delta, \Delta+1\}$, however determining $\chi'(G)$ exactly is known to be NP-hard, even for simple cubic $G$ \cite{H}. A second lower bound for $\chi'(G)$ is given by the density parameter $\lceil\rho(G)\rceil$, where
$$\rho(G)=\max\left\{\tfrac{2|E(G[S])|}{|S|-1}: S\subseteq V(G), |S| \textrm{ odd and } \geq 3 \right\}$$
(to see this notice that $\tfrac{|S|-1}{2}$ is the maximum size of a colour class in $G[S]$).
Goldberg \cite{G-rus} and Seymour \cite{Se} conjectured that in general, $\chi'(G)\in\{\Delta, \Delta+1, \lceil\rho(G)\rceil\}$ (or equivalently, $\chi'(G)\leq \max\{\lceil\rho(G)\rceil, \Delta+1\}$).
If this conjecture is true, then it would imply a polynomial-time algorithm for first checking whether or not $\chi'(G)>\Delta+1$, and then, if the answer is yes, computing $\chi'(G)$ exactly -- in contrast to the NP-hardness of edge-colouring. This implication is due to Edmonds \cite{E}; see also \cite{chap} for an explanation. An asymptotic version of the Goldberg-Seymour conjecture is known (Kahn \cite{Ka}), and the conjecture is known to be true for all graphs not containing a $K_5^-$-minor (Marcotte \cite{Mar}). Most other results towards the conjecture use the \emph{method of Tashkinov trees} to establish an approximation -- namely, to establish a result of the form $\chi'(G)\leq\max \{\lceil\rho\rceil, \Delta+s\}$   for some $s> 1$. For example, Scheide \cite{SchJCTB} proved such an approximation with $s=\frac{\Delta+12}{14}$; see Stiebitz et al. \cite{SSTF} for a complete survey. In this paper we present a specific limit to the method of Tashkinov trees. We also provide a new Tashkinov tree extension.

Suppose that you want to prove $\chi'(G)\leq\max \{\lceil\rho\rceil, \Delta+s\}$ for some $s\geq 1$, and to this end, you assume that $\chi'(G)> \Delta+s$. You succeed if you can find a subgraph $H$ of $G$ that is so dense that it requires more than $\chi'(G)-1$ colours, that is, with $\lceil\rho(H)\rceil >\chi'(G)-1$. This is because such a subgraph would imply that
$$\chi'(G)\geq \lceil\rho(G)\rceil \geq \lceil\rho(H)\rceil >\chi'(G)-1,$$
and hence that $\chi'(G)=\lceil\rho(G)\rceil$, as desired. The method of Tashkinov trees, developed by Tashkinov \cite{T} in 2000 as a common generalization of earlier work by Kierstead \cite{K1} and Vizing \cite{V2}, suggests a candidate for $H$. In particular, starting with any partial edge-colouring of $G$, it describes how to perform (a polynomial number of) \emph{Kempe changes} in order to find a candidate $H$. Given a partial
$k$-edge-colouring of $G$, any pair of colours $a,b \in \{1, 2, \ldots, k\}$ induces a subgraph of $G$ where each component is a path or an even cycle. A \emph{Kempe change} is the action of switching $a$ and $b$ on any such component.

Let $\varphi$ be a partial $k$-edge-colouring of $G$ with $k\geq \Delta+1$, which leaves at least one edge $e_0$ of $G$ uncoloured. Define the set $W_{e_0}^{\varphi}\subseteq V(G)$ recursively as follows:
\begin{enumerate}
\item The two ends of $e_0$ are in $W_{e_0}^{\varphi}$
\item Add a vertex $v$ to $W_{e_0}^{\varphi}$ if there is a colour $\alpha$ such that (when considering $\varphi$): there is an edge $e_v$ joining a vertex in $W_{e_0}^{\varphi}$ to $v$ that is coloured $\alpha$, and there is a vertex $u_v \in W_{e_0}^{\varphi}$ that is missing $\alpha$ (i.e. $\alpha$ is not used on any edge incident to $u_v$).
\end{enumerate}
At any point during the recursive process, the current vertices in $W_{e_0}^{\varphi}$, along with $e_0$ and the edges that have been used in step 2, namely $\{e_v| v\in W_{e_0}^{\varphi}\}$, form a tree called a \emph{$\varphi$-Tashkinov tree starting at $e_0$}. While there may be different maximal $\varphi$-Tashkinov trees starting at $e_0$ (owing to different possible choices for $e_v$), the vertex set of such a maximal tree, namely $W_{e_0}^{\varphi}$, is unique. \emph{Tashkinov's Theorem} \cite{T} says that either $W_{e_0}^{\varphi}$ is $\varphi$-elementary (no pair of vertices in $W_{e_0}^{\varphi}$ have a common missing colour), or that there is a sequence of Kempe changes that will modify $\varphi$ so that both ends of $e_0$ have a common missing colour, and hence the colouring can be extended to $e_0$.

If $W_{e_0}^{\varphi}$ is $\varphi$-elementary, then each colour in $\varphi$ either induces a near-perfect matching in $G[W_{e_0}^{\varphi}]$, or it is \emph{defective}, that is, it occurs on more than one edge between $W_{e_0}^{\varphi}$ and $V(G)\setminus W_{e_0}^{\varphi}$. If there are no defective colours, then according to the definition of $\rho$,
$$\chi'(G)\geq \lceil\rho(G[W_{e_0}^{\varphi}])\rceil \geq \left\lceil\frac{2 \left(k\left(\frac{|W_{e_0}^{\varphi}|-1}{2}\right)+1\right)}{|W_{e_0}^{\varphi}|-1} \right\rceil > k,$$
where the ``$+1$'' comes from the uncoloured edge $e_0$.
In this situation, if $k<\chi'(G)-1$, then one should add a new colour to $\varphi$ and start the analysis again (i.e. consider again the resulting $W_{e_0}^{\varphi}$).
On the other hand, if there are no defective colours and $k=\chi'(G)-1$, then the above sequence of inequalities proves that $\chi'(G)=\lceil\rho(G)\rceil$, as desired.

We now know that in order to prove  $\chi'(G)\leq\max \{\lceil\rho\rceil, \Delta+s\}$ via the method of Tashkinov trees, we ``just'' need to prove that, assuming $\chi'(G)>\Delta+s$,  $W_{e_0}^{\varphi}$ has no defective colours, for $e_0$ some edge left uncoloured by $\varphi$, and $\varphi$ some partial $(\chi'(G)-1)$-edge-colouring with maximum domain (i.e. colouring the most edges among all partial $(\chi'(G)-1)$-edge-colourings of $G$). When $s$ is relatively large, this may not be too difficult: for any such colouring, the total number of colours in $\varphi$ is  $\chi'(G)-1$, however each vertex in $W_{e_0}^{\varphi}$ has at least $\chi'(G)-1-\Delta\geq 1$ missing colours (all of which are distinct since the vertices of $W_{e_0}^{\varphi}$ are $\varphi$-elementary). If there is a defective colour, then it is in addition to the missing colours (by the maximality of our tree), so we get that
$$\chi'(G)-1 \geq |W_{e_0}^{\varphi}|(\chi'(G)-1-\Delta)+2+1,$$
where the ``$+2$'' is from the extra colours missing at the ends of $e_0$. Rearranging this equation, we get
$$\chi'(G)\leq \Delta+ 1 + \frac{\Delta-3}{|W_{e_0}^{\varphi}|-1}.$$
Since $\chi'(G)>\Delta+s$, the larger $|W_{e_0}^{\varphi}|$ and/or $s$ is, the tighter this bound will be. If we are able to argue that $|W_{e_0}^{\varphi}|$ is sufficiently large as to make this inequality invalid, then we get our desired approximation. Here, all of the work is in building a large Tashkinov tree.

As $s$ gets smaller, it becomes more difficult to argue that a sufficiently large Tashkinov tree exists. In particular, it is easy to construct examples for $s=1$ where $W_{e_0}^{\varphi}$ does have defective colours. Here, the natural approach is to do some Kempe changes to modify $\varphi$ so that the new $W_{e_0}^{\varphi}$ is larger and does not have any defective colours, and this seems to work well in practice. However, is such a sequence of Kempe changes always possible? That is, 
can the basic method of Tashkinov trees described above possibly capture the density required to prove the Goldberg-Seymour Conjecture in general? In Section 2 we present an example which answers this question in the negative (Theorem \ref{counterthm}). In fact the second author claimed such an example in her PhD thesis \cite{M}, but a flaw was later found. Here we have a completely new example with no such failing.

The example of Section 2 does not mean the method of Tashkinov trees should be abandoned, but it does tell us that, without a technique beyond Kempe changes for modifying colourings, we cannot just rely on $W_{e_0}^{\varphi}$ to induce our dense subgraph $H$. Other authors, in particular Favrholt, Stiebitz and Toft \cite{FST}, have already done work to define a set $W\supset W_{e_0}^{\varphi}$ that is $\varphi$-elementary and hence provides an improved candidate for $H$. We establish a new such result in Section 3.

\section{A limiting example}\label{limit}

\begin{figure}
\centering
\includegraphics[width=\textwidth]{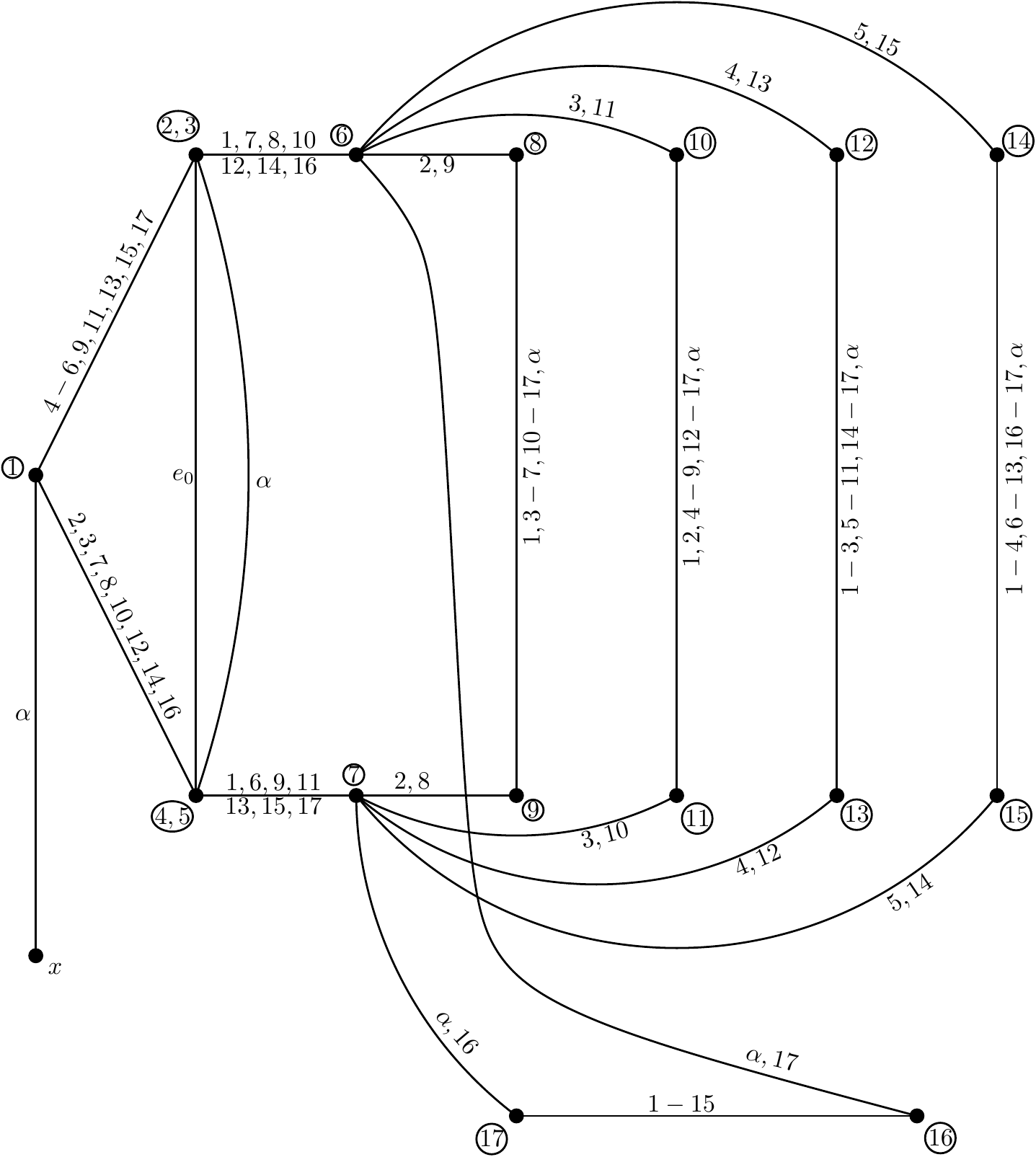}
\caption{A graph and and a colouring where density cannot be captured by a Tashkinov tree, even with an unlimited number of Kempe changes.}
\label{counter}
\end{figure}

The subject of this section is the graph $G$ and partial edge-colouring $\varphi$ given in Figure~\ref{counter}. Every edge in $G$ is coloured except for the edge $e_0$, and in general a single line in the picture represents a number of parallel edges in $G$, whose colours are written along that line.  The encircled colours at the vertices are those colours missing at each vertex. The only vertex without encircled colours is also the only labelled vertex in the graph, $x$, and every colour except for $\alpha$ is missing at $x$. There are 18 colours used in total in this graph, the colours 1-17, and the colour $\alpha$. For each of these 18 colours, let $v_i$ denote the vertex missing colour $i$ in the figure (so $v_2=v_3$ and $v_4=v_5$ by this convention). The following theorem provides a negative answer to the question asked in the introduction, that is, it establishes a firm limit to the method of Tashkinov trees.

\begin{theorem}\label{counterthm} Let $G$, $\varphi$, and $e_0$ be as in Figure~\ref{counter}. Then $\chi'(G)\geq \Delta+2$, and $\varphi$ is a partial $(\chi'-1)$-edge-colouring of $G$ with maximum domain. However, for any $\psi$ obtained from $\varphi$ by a sequence of Kempe changes, $W_{e_0}^{\psi}$ has a defective colour and has $$\lceil\rho(G[W_{e_0}^{\psi}])\rceil<\chi'(G).$$
\end{theorem}

\begin{proof} We can see that $\lceil\rho(G)\rceil=19$ with $V(G)-x$ being the only subset of vertices achieving this density. Since $G$ is 19-colourable (using $\varphi$ and one extra colour for $e_0$), we get that $\chi'(G)=\lceil\rho(G)\rceil=19$. Note that as $\Delta(G)=17$, this means that $\chi'(G)=\Delta+2$ and that $\varphi$ is a $(\chi'(G)-1)$-edge-colouring of $G-{e_0}$.

In Figure~\ref{counter}, $W_{e_0}^{\varphi}=V(G)\setminus\{x, v_{16}, v_{17}\}$, with defective colour $\alpha$. We will show that even after any number of Kempe changes, we cannot increase the size of this set, and in particular, we can never achieve $W_{e_0}^{\psi}=V(G)-x$ for any $\psi$ obtained from $\varphi$ by Kempe changes. (Note that this means that any such $W^{\psi}_{e_0}$ has $\lceil\rho(G[W_{e_0}^{\psi}])\rceil<\chi'(G)$ and hence must have a defective colour.)

We consider the subgraph induced by each pair of colours in $\varphi$. The possible components for such a subgraph are: a copy of $2K_2$ (a doubled edge), an alternating path, or an even cycle of length longer than two.

\begin{claim}\label{notwopaths} No pair of colours in $\varphi$ induce a subgraph containing more than one alternating path.
\end{claim}

\begin{proofc} Suppose some pair of colours did induce two alternating paths. Then the four distinct ends of these paths are each missing one of the two colours in question. However, the colour $\alpha$ is missing only at the vertex $x$, and every other colour is missing at exactly one vertex in $G-x$ and also at $x$. So, these four distinct vertices are not possible.
\end{proofc}

Given Claim \ref{notwopaths}, if a pair of colours does not induce at least one cycle of length longer than two, then any Kempe change involving those colours is merely a relabelling of $\varphi$. Most pairs of colours in $\varphi$ fall into this category. To see this, it is useful to use the symmetry of $\varphi$ and group some of the colours together. We let $A=\{6, 7\}$, $B=\{16,17\}$, $C=\{8, 9, \ldots, 15\}$, and $U=\{2, 3, 4, 5\}$, but treat $1$ and $\alpha$ on their own. Let $a, a'\in A$, $b, b'\in B$, $c, c'\in C$ and $u, u' \in U$. Then a manual check verifies that the subgraph induced by any of the following pairs of colours consists of a single alternating path and copies of $2K_2$: $\{a, a'\}$, $\{a,b\}$, $\{a,c\}$, $\{a,u\}$, $\{a,1\}$, $\{a, \alpha\}$, $\{b, b'\}$, $\{b,c\}$, $\{b,u\}$, $\{b,1\}$, $\{b, \alpha\}$, $\{c,c'\}$, $\{c,u\}$, $\{c,1\}$, $\{c, \alpha\}$, $\{u,1\}$.

It remains only to consider the pairs $\{u, u'\}$, $\{u, \alpha\}$, $\{1, \alpha\}$, each of which does induce a cycle of length longer than 2. The easiest of these cases is $\{u, u'\}$  equal to $\{2,3\}$ or $\{4,5\}$. In this situation, the induced subgraph has no alternating path and just consists of a single alternating cycle and copies of $2K_2$. So again here, any Kempe change is just a relabelling of $\varphi$. In the remaining cases, the subgraph induced contains both an alternating path and an alternating cycle of length longer than two, so there are Kempe changes which are not just relabellings of $\varphi$.

Consider first $\{u, u'\}$ where $u\in \{2, 3\}$ and $u'\in\{4,5\}$. The subgraph induced by these colours
contains the path $v_2, v_1, v_4$, copies of $2K_2$, and a 6-cycle containing the pairs of vertices $v_6, v_7$, $v_i, v_{i+1}$, and $v_j, v_{j+1}$ for some distinct $i,j\in\{8, 10, 12, 14\}$. Note that we get an isomorphic copy of $G$, say $G'$, by mapping $v_{i}$ and $v_{j}$  to each other and mapping $v_{i+1}$ to $v_{j+1}$ to each other. We will argue that doing the $(u, u')$-switch along the 6-cycle (i.e. switching the colours $u$ and $u'$ along the 6-cycle) results in a colouring which is not a relabelling of $\varphi$ on $G$, but is in fact a relabelling of $\varphi$ on the isomorphic graph $G'$.

To this end, consider $\varphi$ on $G$, and imagine first executing the $(i,j)$-switch along the alternating path starting at $v_{i}$  (which ends at $v_j$) and then the $(i+1, j+1)$-switch along the alternating path starting at $v_{i+1}$ (which ends at $v_{j+1}$). We have already concluded that both of these switches result in merely a relabelling of $\varphi$ on $G$, so the result of both actions, say $\varphi_0$, is also a relabelling of $\varphi$ on $G$. Now apply the isomorphism described above to get $\varphi_0$ on the isomorphic graph $G'$. This is exactly the same result as if we did the $(u, u')$-switch along the 6-cycle to modify $\varphi$ in $G$.

Consider now the pair $\{u, \alpha\}$. The subgraph induced by this pair contains a path with vertex set $\{x, v_1, v_2, v_4\}$, copies of $2K_2$, and a 6-cycle containing the pairs of vertices $v_6, v_7$, $v_{i}, v_{i+1}$, and $v_{16}, v_{17}$ for some $i\in\{8, 10, 12, 14\}$. We can make an argument identical to that of the previous case, except that instead of $v_j$ and $v_{j+1}$ (and the colours $j$ and $j+1$) we have instead the vertices $v_{16}$ and $v_{17}$ (and the colours 16 and 17), respectively.

It remains now only to consider the subgraph induced by $\{1, \alpha\}$. The induced subgraph consists of the single-edge path from $x$ to $v_1$, copies of $2K_2$, and the 6-cycle $v_2, v_4, v_7, v_{17}, v_{16}, v_6$. By doing a Kempe change on either the cycle or the path, we actually do get a new colouring (say $\varphi'$), which is not a relabelling of $\varphi$ in $G$ or in any isomorphic copy of $G$. We thus proceed with caution.

Without loss of generality, let us say that $\varphi'$ has been obtained by switching $1$ and $\alpha$ along the path $xv_1$. Then $W_{e_0}^{\varphi'}=\{v_1, v_2, v_4\}$, which is obviously too small to capture the required density of 19. However, we must also consider any sequence of Kempe changes that could be done to $\varphi'$ which might allow this set to increase in size. If such a sequence exists, we must be able to start with a switch involving exactly one of $1$ or $\alpha$. We will show that every such switch is either just a relabelling of $\varphi'$ on $G$, or a relabelling of $\varphi'$ on an isomorphic copy of $G$, and this will complete our proof.

The pairs $\{a, \alpha\}$, $\{b,\alpha\}$ and $\{c, \alpha\}$ induce, in $\varphi$, a single path ending at $x$ along with some copies of $2K_2$. The subgraph they induce in $\varphi'$ is different only in that each path is truncated at $v_1$. Similarly, the pairs $\{a,1\}$, $\{b, 1\}$,  $\{c, 1\}$, $\{u, 1\}$ induce, in $\varphi$, a single path ending at $v_1$ along with some copies of $2K_2$. The subgraph they induce in $\varphi'$ is different only in that the path is extended to $x$ and there are no isolated vertices. So, in all of these cases, we get only a relabelling of $\varphi'$ by making a switch.

The only remaining pair to consider is $\{u, \alpha\}$. Recall that in $\varphi$, the subgraph induced by this pair contains a path with vertex set $\{x, v_1, v_2, v_4\}$, copies of $2K_2$, and a 6-cycle containing the pairs of vertices $v_6, v_7$, $v_{i}, v_{i+1}$, and $v_{16}, v_{17}$ for some $i\in\{8, 10, 12, 14\}$.
The difference in $\varphi'$ is just that the path ends at $v_1$  instead of $x$. Our above argument about switching on the $\{u, \alpha\}$ cycle in $\varphi$ therefore works identically in $\varphi'$. So, we conclude that such a switch would just be a relabelling of $\varphi'$ on a graph isomorphic to $G$.
\end{proof}

In the above example of Figure \ref{counter}, it may seem as though the Kempe changes were insufficient because of the placement of $e_0$. Suppose that, in addition to Kempe changes, we allow $e_0$ to be moved -- that is, for the colour from some edge to be removed provided it allows for $e_0$ to be coloured. In $\varphi$, the $(u, u')$-alternating path where $u\in\{2, 3\}$ and $u'\in\{4,5\}$ forms a cycle with $e_0$, and hence $e_0$ could be moved along this cycle (triangle). However, changing $\varphi$ and $e_0$ in this way would have no effect on $W^{\varphi}_{e_0}$. Moreover, without first doing one or more Kempe changes, the only colours we could use to colour $e_0$ would be one of $\{2, 3, 4, 5\}$, and the effect would be the same as moving the edge around the triangle $v_1, v_2, v_4$. If we hope to first do some Kempe changes to facilitate a different move for $e_0$, we have only to look in the above proof to see that this would not make any difference.

\section{A larger elementary set}

Let $G$ be a loopless graph with $\chi'(G)\geq \Delta+2$, let $\varphi$ be a partial $k$-edge-colouring of $G$ with $k\geq \Delta+1$ colours, and let $e_0$ be an edge uncoloured by $\varphi$. Consider $W^{\varphi}_{e_0}$  and define $\mathcal{M}^{\varphi}_{e_0}$ to be the set of all colours missing at (at least one) vertex in $W^{\varphi}_{e_0}$. Also define $\mathcal{D}^{\varphi}_{e_0}$ to be the set of defective colours for $W^{\varphi}_{e_0}$. Note that $\mathcal{D}^{\varphi}_{e_0}\cap \mathcal{M}^{\varphi}_{e_0}=\emptyset$ by definition of $W^{\varphi}_{e_0}$. Given $W^{\varphi}_{e_0}$, fix $T$ to be any maximal $\varphi$-Tashkinov tree starting at $e_0$, and define $\mathcal{U}^{\varphi}_{T}$ to be the set of colours used (i.e. appearing) on the edges of $T$.

The following result, which can be seen as an algorithmic approach to colouring via the method of Tashkinov trees, contains: (a) Tashkinov's Theorem \cite{T}; (b) work by Favrholt, Steibitz and Toft \cite{FST} (see also Theorem 5.11 of \cite{SSTF}), and; (c) a new contribution which extends (b). It should be mentioned that Favrholt, Stiebitz and Toft have a different extension to (b) (via a fan structure and the concept of absorbing vertices) which is described in Theorem 5.18 of \cite{SSTF}.

\begin{figure}
\centering
\includegraphics[width=.55\textwidth]{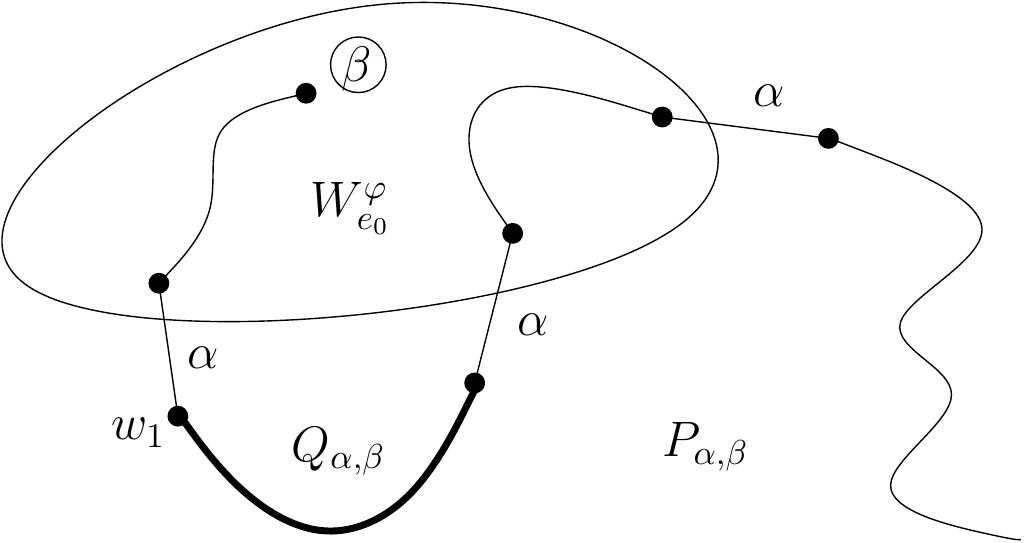}
\caption{The paths $P_{\alpha, \beta}$ and $Q_{\alpha, \beta}$.}
\label{PandQ}
\end{figure}

\begin{theorem}\label{collection} Let $G$ be a loopless graph with $\chi'(G)\geq \Delta+2$, let $\varphi$ be a partial $k$-edge-colouring of $G$ with $k\geq \Delta+1$ colours, and let $e_0$ be an edge uncoloured by $\varphi$.
If none of the following statements are true:
\begin{enumerate}
\item[(1)] There is a sequence of Kempe changes that transform $\varphi$ into a colouring $\varphi'$ where the two ends of $e_0$ have a common missing colour (and hence $\varphi'$ can be extended to $e_0$).
\item[(2)] There is a sequence of Kempe changes that transform $\varphi$ into a colouring $\varphi'$ where $|W^{\varphi'}_{e_0}|>|W^{\varphi}_{e_0}|$.
\item[(3)] $\chi'(G)\geq \lceil\rho(G[W^{\varphi}_{e_0}])\rceil > k$.
\end{enumerate}
then all of the following statements are true:
\begin{enumerate}
	\item[$(a)$] $W^{\varphi}_{e_0}$ is $\varphi$-elementary and has at least one defective colour, with each defective colour appearing an odd number of times between $W^{\varphi}_{e_0}$ and $G\setminus W^{\varphi}_{e_0}$.
	\item[$(b)$] If $\beta \in \mathcal{M}^{\varphi}_{e_0}\setminus \mathcal{U}^{\varphi}_T$ and $\alpha \in \mathcal{D}^{\varphi}_{e_0}$, the maximal $(\alpha, \beta)$-alternating path $P_{\alpha, \beta}$ that starts in $W^{\varphi}_{e_0}$ must contain all $\alpha$-edges having exactly one end in $W^{\varphi}_{e_0}$. Moreover, if $w_1, w_2$ are the first two vertices of $P_{\alpha, \beta}$ that are not in $W^{\varphi}_{e_0}$, then $W^{\varphi}_{e_0}\cup \{w_1, w_2\}$ is $\varphi$-elementary.
	\item [(c)] Let $Q_{\alpha, \beta}$ be the maximal-length segment of $P_{\alpha, \beta}$ which starts at $w_1$ and contains no vertices of $W^{\varphi}_{e_0}$. Then $W^{\varphi}_{e_0} \cup Q_{\alpha,\beta}$ is $\varphi$-elementary. See Figure \ref{PandQ}.
\end{enumerate}
\end{theorem}

\begin{proof} Suppose that none of the statements in the first list are true. Then, since (1) is not true, Tashkinov's Theorem tells us that $W_{e_0}^{\varphi}$ is $\varphi$-elementary. In particular, this means that $W_{e_0}^{\varphi}$ has an odd number of vertices, so any defective colour must appear on an odd number of edges between $W_{e_0}^{\varphi}$ and $G\setminus W_{e_0}^{\varphi}$. Since $W_{e_0}^{\varphi}$ is $\varphi$-elementary and (3) is not true, our argument from the introduction tells us that $W_{e_0}^{\varphi}$ must have at least one defective colour. Hence (a) is established. Favrholt, Stiebitz and Toft \cite{FST}(see also Theorem 5.11 of \cite{SSTF}) established (b), for which the assumption that (2) is not true is needed. We now extend this by proving (c). To this end, the following claim is useful.

\begin{claim} $|\mathcal{M}^{\varphi}_{e_0}\setminus \mathcal{U}^{\varphi}_T| \geq 4$.
\end{claim}
\begin{proofc} We know that $$|\mathcal{M}^{\varphi}_{e_0}| \geq |W_{e_0}^{\varphi}| (k-\Delta) + 2 \geq |W_{e_0}^{\varphi}|+2,$$
where the $+2$ comes from the fact that $e_0$ is uncoloured. We also know that $$|\mathcal{U}^{\varphi}_T| \leq |E(T)|-1=|W_{e_0}^{\varphi}|-2.$$ Hence we have our desired result.
\end{proofc}

Let the vertices of $Q_{\alpha, \beta}$ be $w_1, w_1, \ldots, w_{q}$ (note $w_1, w_2$ are as in (b)). Suppose first that there exists $\varepsilon \in \overline{\varphi}(w_i) \cap \overline{\varphi}(v)$ for some $v\in W_{e_0}^{\varphi}$ and some $i\in\{1, 2, \ldots, q\}$. We show by induction on $i$ that this cannot occur. If $i=1$ or $2$ then (b) tells us this immediately. Now suppose that $i\geq 3$. We may assume, without loss of generality, that $\varepsilon\not\in\mathcal{U}^{\varphi}_T$, as if $\varepsilon\in\mathcal{U}^{\varphi}_T$ we can switch the names of $\varepsilon$ and some colour in $\mathcal{M}^{\varphi}_{e_0}\setminus(\mathcal{U}^{\varphi}_T\setminus\{\beta\})$ (which is not an empty set by the above claim) everywhere in $G[W_{e_0}^{\varphi}]$. Now choose $\delta\in\overline{\varphi}(w_{i-1})$ (that is, choose a colour $\delta$ that is missing at $w_{i-1}$). Since $\varphi$ has at least $\Delta+1$ colours, we can do this. Note that, by induction, $\delta\not\in\mathcal{M}^{\varphi}_{e_0}$. Consider the maximal $(\varepsilon, \delta)$-alternating path $L$ beginning at $w_{i-1}$, and define $\varphi'$ to the the edge-colouring obtained by switching $\varepsilon$ and $\delta$ along $L$. Note that this does not change the colours of any edges in $T$ (since neither colour was in $\mathcal{U}^{\varphi}_T$) nor $P_{\alpha, \beta}$, nor does it change the colours missing at any vertices in $W_{e_0}^{\varphi}$, aside from perhaps $v$ (if $L$ ends at $v$). If $L$ does not end at $v$, then after the switch we have our desired result by induction since $w_{i-1}$ and $v$ are now both missing $\varepsilon$. So, suppose that $L$ does end at $v$. In this case, after the switch, $\varepsilon \in \overline{\varphi}(w_{i-1})\cap \overline{\varphi}(w_i)$. By recolouring the edge between these two vertices with $\varepsilon$, $P_{\alpha, \beta}$ will end at $w_{i-1}$, violating (b).

Suppose now that there exists $\varepsilon \in \overline{\varphi}(w_i) \cap \overline{\varphi}(w_j)$ for some $i, j\in\{1, 2, \ldots, q\}$ with $i<j$. We show that this cannot occur with an induction on $j$, and a secondary induction on $j-i$. Note that by the result we have just established, we know that $\varepsilon\not\in\mathcal{M}^{\varphi}_{e_0}$. By (b) we may assume that $j\geq 3$. If $j-i=1$, then we do as we did at the end of the previous case: we recolour the edge between the two vertices, truncating $P_{\alpha, \beta}$ for a contradiction to (b). Hence we may assume that $j-1\geq 2$. Choose $\delta\in\overline{\varphi}(w_{i+1})$. By our above work we know that $\delta\not\in\mathcal{M}^{\varphi}_{e_0}$. Consider the maximal $(\varepsilon, \delta)$-alternating path $L$ starting at $w_{i+1}$. Define $\varphi'$ by switching $\varepsilon$ and $\delta$ along $L$. Note that since both $\varepsilon$ and $\delta$ are not in $M^{\varphi}_{e_0}$ (nor are they $\alpha$ or $\beta$),  this switch does not change the colours of any edges in $T$ or $P_{\alpha, \beta}$, nor does it change the colours missing at any vertices in $W_{e_0}^{\varphi}$. After the switch, we either complete our proof by primary induction hypothesis with $w_i$ and $w_{i+1}$ (if $L$ does not end at $w_i$), or we complete our proof by secondary induction hypothesis with $w_{i+1}$ and $w_j$ (if $L$ ends at $w_i$).
\end{proof}

\begin{figure}
\centering
\includegraphics[width=\textwidth]{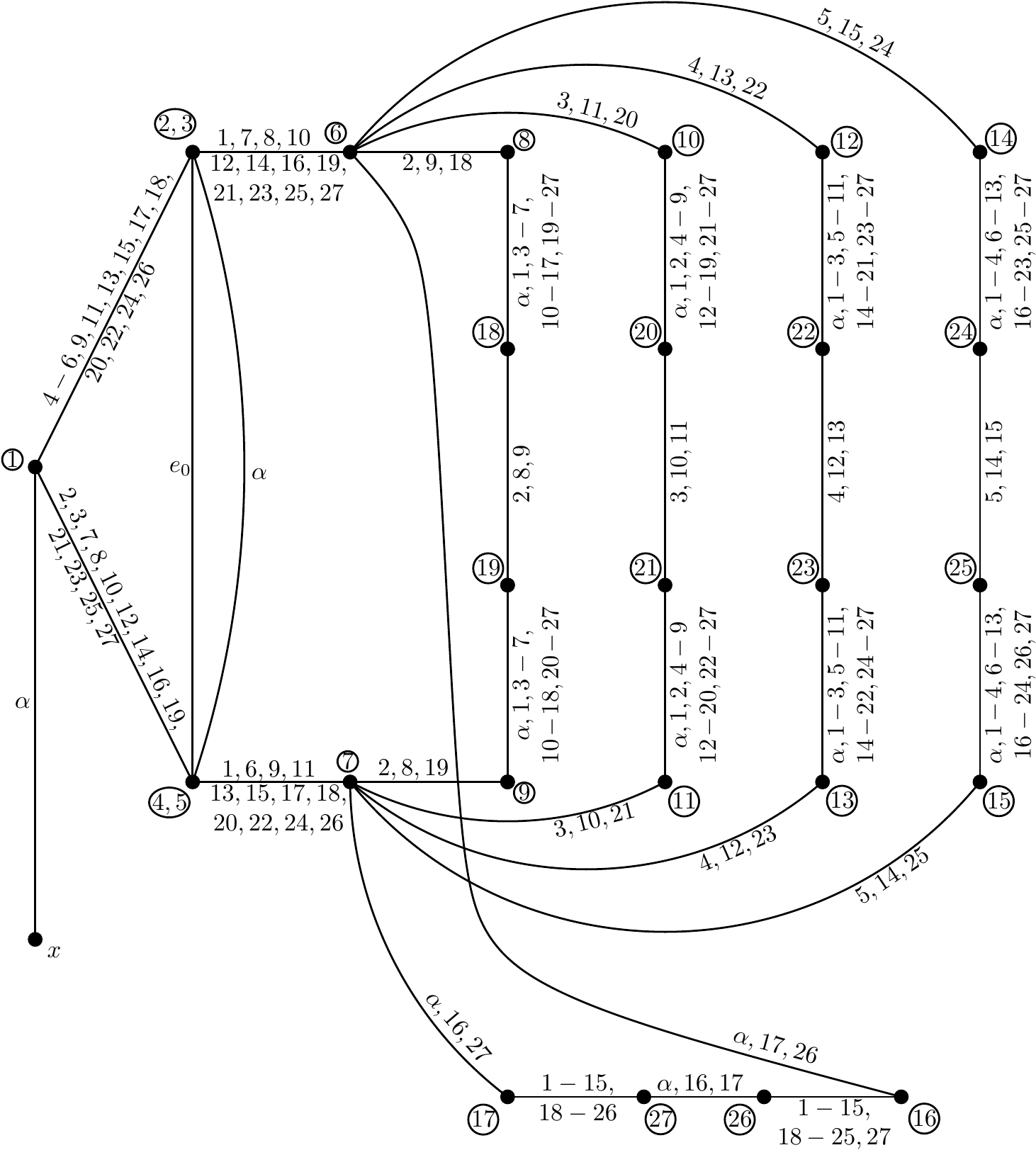}
\caption{A graph and and a colouring where density cannot be captured by the structure in Theorem~\ref{collection}(c), even with an unlimited number of Kempe changes.}
\label{counterexample2}
\end{figure}

Theorem \ref{counterthm} proves that, for the example in Figure \ref{counter}, density cannot be captured by a Tashkinov tree, even if we allow for an unlimited number of Kempe changes. However, the structure in Theorem ~\ref{collection}(b) \emph{does} capture density for that example. To see this, note that there, $\mathcal{M}^{\varphi}_{e_0}=\{1, 2, \ldots, 15\}$, $\mathcal{D}^{\varphi}_{e_0}=\{\alpha\}$, and we can take $T$ so that $\mathcal{U}^{\varphi}_{T}=\{1, 2, 3, 4, 5\}$. So by choosing $\beta=6$, for example, we get that $P_{\alpha, \beta}$ has $w_1=v_{16}$ and $w_2=v_{17}$. Hence part (b) tells us that
$W^{\varphi}_{e_0}\cup \{w_1, w_2\}=V(G)\setminus x$ is $\varphi$-elementary. This set does capture density since $\lceil \rho(V(G)\setminus x) \rceil =\chi'(G)$.  We can, however, modify this example (to Figure \ref{counterexample2}) so that (b) is no longer sufficient, although (c) is. In Figure \ref{counterexample2}, we use the same conventions as in Figure \ref{counter} in terms of indicating the colours used on each edge and missing at each vertex, and in labelling a vertex $v_i$ if the colour $i$ is missed there.

\begin{theorem}\label{g2} Let $G$,$\varphi$, $e_0$ be as in Figure~\ref{counterexample2}. Then $\chi'(G)\geq \Delta+2$, and $\varphi$ is a partial $(\chi'-1)$-edge-colouring of $G$ with maximum domain. Then, for any
 $\alpha, \beta$ chosen as in Theorem \ref{collection}, $\chi'(G)=\lceil\rho(G[W_{e_0}^{\varphi}\cup Q_{\alpha,\beta}])\rceil$. On the other hand, for any $\psi$ obtained from $\varphi$ by a sequence of Kempe changes, $W_{e_0}^{\psi}$ has a defective colour and has  $$\lceil\rho(G[W_{e_0}^{\psi}\cup\{w_1, w_2\}])\rceil< \chi'(G),$$
for any $w_1, w_2$ chosen as in Theorem \ref{collection} with respect to $\psi$.
\end{theorem}

\begin{proof} We can see that $\lceil\rho(G)\rceil=29$ with $V(G)-x$ being the only subset of vertices achieving this density. Since $G$ is 29-colourable (using $\varphi$ and one extra colour for $e_0$), we get that $\chi'(G)=\lceil\rho(G)\rceil=29$. Note that as $\Delta(G)=27$, this means that $\chi'(G)=\Delta+2$ and that $\varphi$ is a $(\chi'(G)-1)$-edge-colouring of $G-{e_0}$.

In Figure~\ref{counterexample2}, $W_{e_0}^{\varphi}=V(G)\setminus\{x, v_{16}, v_{17}, v_{26}, v_{27}\}$,
$\mathcal{M}^{\varphi}_{e_0}=\{1, 2, \ldots, 15,18,19,\ldots,25\}$, and $\mathcal{D}^{\varphi}_{e_0}=\{\alpha\}$. To form a $\varphi$-Tashkinov tree $T$ starting at $e_0$ with vertex set $W_{e_0}^{\varphi}$, we need to use all of the colours 1, 2, 3, 4, 5 on the edges of $T$, i.e., $\{1, 2, 3, 4, 5\}\subseteq \mathcal{U}_T^{\varphi}$. Since we need not use any other colours to build such a $T$, we can choose for $\beta$ any colour in $\{6, 7, \ldots, 15,18,19,\ldots,25\}$. For any such choice $Q_{\alpha, \beta}$ is the same, namely $Q_{\alpha, \beta}=(v_{17}, v_{27}, v_{26}, v_{16})$. Hence for any such choice of $\beta$, $\lceil\rho(G[W_{e_0}^{\psi}\cup Q_{\alpha,\beta}])\rceil=\chi'(G)$. On the other hand, since $|W_{e_0}^{\varphi}|=21$, for any such choice $|W_{e_0}^{\varphi}\cup\{w_1, w_2\}|=23<25=|V(G)-x|$ and hence $\lceil\rho(G[W_{e_0}^{\varphi}\cup\{w_1, w_2\}])\rceil< \chi'(G)$. We will show that for any $\psi$ obtained from $\varphi$ by a sequence of Kempe changes, $|W_{e_0}^{\psi}|\leq 21$.

It is useful to describe how we made Figure \ref{counterexample2} from Figure \ref{counter}. First, for each edge of colour $i\in\{8, 10, 12, 14, 16\}$ (in Figure \ref{counter}) we made a parallel edge and gave it colour $i+11$. Then, for each edge of colour $i+1\in\{9, 11, 13, 15, 17\}$ we made a parallel edge and gave it colour $i+10$. This process yielded the new colours $18, \ldots, 27$. Then, we took the set of the edges $E_i$ between $v_i$ and $v_{i+1}$ (for $i\in\{8,10,12,14,16\}$) and made two copies of $E_i$: we deleted the original set, but placed one copy of $E_i$ between $v_i$ and the new vertex $v_{i+10}$, and place the other copy of $E_i$ between $v_{i+1}$ and $v_{i+11}$; to the first copy we added a parallel edge coloured $i+11$, and to the second copy we added a parallel edge coloured $i+10$. Finally we put three edges between $v_i$ and $v_{i+1}$ (for $i\in\{8,10,12,14,16\}$), and coloured these three edges with the colours $i$, $i+1$, and  $\tfrac{i}{2}-2$ , except in the case of $i=16$, where we used the colour $\alpha$ in place of $\tfrac{i}{2}-2$.

Given the way that we made Figure \ref{counterexample2} from Figure \ref{counter}, every pair of colours in $\{\alpha, 1, 2, \ldots, 17\}$ (i.e. every pair of colours from Figure \ref{counter}) induce as components the same number of paths and the same number of long (length more than 2) cycles in Figure \ref{counterexample2} as they did in Figure \ref{counter}, although these paths and even cycles may be longer than they were before. Hence, given the argument in the proof of Theorem \ref{counterthm}, the only such pairs inducing both a path component and a long cycle component in Figure \ref{counterexample2} are $\{u, u'\}, \{u, \alpha\}, \{u', \alpha\}$, and  $\{1, \alpha\}$, where $u\in\{2,3\}$ and $u'\in\{4,5\}$.

In Figure \ref{counterexample2}, view each colour $i\in\{8, 10, 12, 14, 16\}$ as being twinned with the colour $i+11$, and each colour $i+1\in\{9, 11, 13, 15, 17\}$ as twinned with the colour $i+10$. Each such twinning induces only copies of $2K_2$, except for one path of length 2. Because of this, a pair of colours from $\{18, \ldots, 27\}$ induces a subgraph which is the same as the subgraph induced by their twins in $\{8, \ldots, 17\}$, except that the former will have two more copies of $2K_2$ and the path in the latter will be four edges longer. Moreover, given a pair of colours $\{w,y\}$ with $w\in\{\alpha, 1, 2, \ldots, 17\}$ and $y\in\{18, \ldots, 27\}$ (where $w$ and $y$ are not twins), the subgraph induced by $\{w,x\}$ is the same as the subgraph induced by $\{w, y'\}$ (where $y'$ is the twin of $y$), except that the length of the path component and the number of $2K_2$ components will be different.

It remains now only to consider the pairs $\{u, u'\}, \{u, \alpha\}, \{u', \alpha\}$, and  $\{1, \alpha\}$, where $u\in\{2,3\}$ and $u'\in\{4,5\}$. All of these pairs do induce subgraphs consisting of both a path and a long even cycle, in addition to copies of $2K_2$.
In the proof of Theorem \ref{counterthm} we dealt with all of these pairs, save $\{1, \alpha\}$, in the same way: by showing that a switch along the long cycle is equivalent to a relabelling of an isomorphic copy of $G$. We will do this again here, although we need to do slightly more work to get the relabelling.

Consider a pair $\{u, u'\}, \{u, \alpha\}$, or $\{u', \alpha\}$, where $u\in\{2,3\}$ and $u'\in\{4,5\}$. The long cycle induced by this pair is the 10-cycle $(v_6, v_{i}, v_{i+10}, v_{i+11}, v_{i+1}, v_7, v_{j+1}, v_{j+11},\\ v_{j+10}, v_j)$, for some distinct $i,j\in\{8, 10, 12, 14, 16\}$. Note that we get an isomorphic copy of $G$, say $G'$, by mapping $v_{i}, v_{i+10}, v_{i+11}, v_{i+1}$ to $v_{j}, v_{j+10}, v_{j+11}, v_{j+1}$, and vice versa. We will argue that swapping the pair of colours in question around the 10-cycle results in a colouring which is not a relabelling of $\varphi$ on $G$, but is in fact a relabelling of $\varphi$ on the isomorphic graph $G'$.

To this end, consider $\varphi$ on $G$, and imagine executing all of the following switches, in sequence: the $(i,j)$-switch starting at $v_i$ (and ending at $v_j$), the $(i+1, j+1)$-switch starting at $v_{i+1}$ (and ending at $v_{j+1}$), the $(i+10, j+10)$-switch starting at $v_{i+10}$ (and ending at $v_{j+10}$), and the $(i+11, j+11)$-switch starting at $v_{i+11}$ (and ending at $v_{j+11}$). We have already concluded that all of these switches result in merely a relabelling of $\varphi$ on $G$, so the result of all 4 actions, say $\varphi_0$, is also a relabelling of $\varphi$ on $G$. Now apply the isomorphism described above to get $\varphi_0$ on the isomorphic graph $G'$. This is exactly the same result as if we switched along the 10-cycle in question.

It remains now only to consider the pair $\{1, \alpha\}$. As in the proof of Theorem~\ref{counterthm}, by performing a Kempe change on either the path or cycle induced by $\{1,\alpha\}$, we actually get a new colouring (say $\varphi'$) which is not a relabelling of $\varphi$ in $G$ or in any isomorphic copy of $G$. Without loss of generality, let us say $\varphi'$ has been obtained by switching $1$ and $\alpha$ along the path $xv_1$. Then $W_{e_0}^{\varphi'}=\{v_1,v_2,v_4\}$, which is obviously too small to capture the density of $29$. However, we must also consider any sequence of Kempe changes that could be done to $\varphi'$ which might allow this set to increase in size. If such a sequence exists, we must be able to start with a switch involving exactly one of $1$ or $\alpha$. We will show that each such switch is just a relabelling of $\varphi'$ on $G$, or a relabelling of $\varphi'$ on an isomorphic copy of $G$, and this will complete our proof.

The pairs $\{s,\alpha\}$ with $s\in\{6,7,\ldots,27\}$ induce, in $\varphi$, a single path ending at $x$ along with some copies of $2K_2$. The subgraph they induce in $\varphi'$ is different only in that each path is truncated at $v_1$. Similarly, the pairs $\{s,1\}$ induce, in $\varphi$, a single path ending at $v_1$ along with some copies of $2K_2$. The subgraph they induce in $\varphi'$ is different only in that the path is extended to $x$ and there are no isolated vertices. So, in all of these cases, we get only a relabelling of $\varphi'$ by making a switch.

The only remaining pair to consider is $\{x,\alpha\}$ where $y\in\{2, 3, 4, 5\}$. Recall that in $\varphi$, the subgraph induced by this pair contains a path with vertex set $\{y,v_1,v_2,v_4\}$, copies of $2K_2$, and a $10$-cycle $(v_6,v_i,v_{i+10},v_{i+11},v_{i+1},v_7,v_{17},v_{27},v_{26},v_{16})$ for some $i\in\{8,10,12,14\}$. The difference in $\varphi'$ is just that the path ends at $v_1$ instead of $y$. Our above argument about switching on the $\{y,\alpha\}$ cycle in $\varphi$ therefore works identically in $\varphi'$. So, we conclude that such a switch would just be a relabelling of $\varphi'$ on a graph isomorphic to $G$.
\end{proof}

\end{document}